\documentclass[11pt]{amsart}
\usepackage{amsmath}
\usepackage{amssymb}
\usepackage{latexsym}
\usepackage{graphicx}

\newtheorem{theorem}{Theorem}
\newtheorem{corollary}{Corollary}
\theoremstyle{remark}
\newtheorem{example}{Example}

\begin{document}

\title[An extension of Janowski functions with applications]{\large An extension of Janowski functions with applications\footnote{\bf In honour of Professor T. O. Opoola.}}

\author[K. O. Babalola]{K. O. BABALOLA}

\begin{abstract}
Janowski functions are regular functions of the form $p(z)=1+p_1z+\cdots$, which are subordinate to $L_0(a,b:z)=(1+az)/(1+bz)$ for real numbers $a,b$ such that $-1\leq b<a\leq 1$. In this paper, we expand the space of the constants, $a$ and $b$, of the superordinate function, $L_0(a,b:z)$, to a larger set of complex numbers and provide some applications.
\end{abstract}



\maketitle

\section{Introduction}
The greatly generalized family of Janowski functions consists of functions $p(z)=1+p_1z+\cdots$, which are regular in the unit disk $E=\{z\in \mathbb{C}\colon |z|<1\}$ and are subordinate to $L_0(a,b:z)=(1+az)/(1+bz)$ for real numbers $a,b$ such that $-1\leq b<a\leq 1$ \cite{JW}. The important roles of the classical Janowski functions in the theory of regular schlicht functions are well known. Practically most classes of geometric functions have representation in terms of them, with the various choices of the parameters, $a$ and $b$, of the superordinate (dominant) function, $$L_0(a,b:z)=\frac{1+az}{1+bz},\;-1\leq b<a\leq 1.$$ See \cite{AA,BA,FD,AR,FA,SA,LA,RM,GI,HM,DL,LM,LX} for example for many works associated with the Janowski functions.\vskip 2mm

It is however well known that the parameters $a$ and $b$ of the Janowski maps have been restricted to the set of real numbers only as defined above. The object of this paper, therefore, is to extend the space of the parameters $a$ and $b$ to a well-defined larger set of complex numbers. The extension of these parameters to the set of complex numbers has provided more possibilities and insights to the study of functions which are analytic and univalent in the unit disk. New relationships between the extended Janowski functions and other classes of functions yielding many new classes of applicable analytic and univalent functions have been identified.\vskip 2mm

It is also important to remark that these new relationships shall open a floodgate to many researchable new problems in the theory of regular schlicht functions. Many future researches will spring up from the new extension of the parameters. We shall demonstrate this with some new definition of classes of functions.\vskip 2mm

We begin by carefully defining a non-vacuous concept as follows:\vskip 2mm

Let $s$ and $t$ be complex numbers such that:\vskip 2mm

{\rm(i)} $|s|\leq 1$, $|t|\leq 1$ and

{\rm(ii)} $\text{Re}\;t<\text{Re}\;s$.\vskip 2mm

We affirm that the concept is non-vacuous as it is not difficult to see from the following example:\vskip 2mm
 
\begin{example} For real numbers $a,b,\lambda$ such that $-1\leq b<a\leq 1$ and $\lambda\in(-1,1)$, define $$s=\frac{a+i\lambda b}{1+i\lambda},\;t=\frac{b+i\lambda a}{1+i\lambda}\eqno{(1.1)}.$$
\end{example}

\begin{proof}
By simple computations, given the conditions on the parameters $a,b$ and $\lambda$, we can see easily that $s$ and $t$ satisfy $|s|\leq 1$, $|t|\leq 1$. Similarly we find that 
$$\text{Re}\;s=\frac{a+\lambda^2 b}{1+\lambda^2},\;\text{Re}\;t=\frac{b+\lambda^2 a}{1+\lambda^2}$$ so that 
$$\text{Re}\;s-\text{Re}\;t=\frac{(a-b)(1-\lambda^2)}{1+\lambda^2}>0$$ which shows that $\text{Re}\;t<\text{Re}\;s$ as desired.  
\end{proof}

The couple (1.1) is a generating string for the complex numbers $s$ and $t$.\vskip 2mm

Now let $\Omega$ be the class of functions $\omega(z)=b_1z+b_2z^2+\cdots$ which are also regular $E$ and satisfy $\omega(0)=0$, $|\omega(z)|<1$ for $z\in E$. Then for complex numbers $s$ and $t$ as defined in (1.1), we say:\vskip 2mm

{\sc Definition.} $P[s,t]$ is the family of functions, $p(z)$, normalized by $p(0)=1$, regular in $E$ and such that for some $\omega\in\Omega$, $$p(z)=\frac{1+s\omega(z)}{1+t\omega(z)},\; z\in E.\eqno{(1.2)}$$

Equivalently $P[s,t]$ is the family of functions $$p(z)=1+p_1z+p_2z^2+\cdots$$ such that $$p(z)\prec L_0(s,t:z)=\frac{1+sz}{1+tz},\;z\in E$$ where $\prec$ denotes subordination.\vskip 2mm

\begin{theorem} With the complex numbers $s$ and $t$ as defined in (1.1), $$L_0(s,t:z)=\frac{1+sz}{1+tz}$$ is well defined (normalized, has positive real part and maps the unit disk $E$ onto a convex domain).
\end{theorem}

\begin{proof}
For some real numbers $a,b$ such that $-1\leq b<a\leq 1$, it is well know that $L_0(a,b:z)=(1+az)/(1+bz)$ is normalized, has positive real part and maps the unit disk $E$ onto a convex domain. Then it is sufficient to show that $L_0(s(\lambda),t(\lambda):z)$ is a linear transformation of the half-plane defined by $L_0(a,b:z)$ onto itself \cite{AW}. Since $s(\lambda)$ and $t(\lambda)$ satisfy the required conditions and so are well defined, then we have $$\aligned L_0(s,t:z)&=\frac{\frac{1+az}{1+bz}+i\lambda}{1+i\lambda\frac{1+az}{1+bz}}\\ &=\frac{L_0(a,b:z)+i\lambda}{1+i\lambda L_0(a,b:z)}.\endaligned$$
which is the linear transformation of the half-plane defined by $L_0(a,b:z)$ onto itself as required.
\end{proof}

Since $L_0(s,t:z)$ is a linear transformation of the half-plane defined by $L_0(a,b:z)$ onto itself, an important corollary to the above result is:\vskip 2mm

\begin{corollary} 
Functions in $P[s,t]$ have positive real parts.
\end{corollary}

The classes $P[s,t]$ corresponds to the well-known Janowski classes of functions, $P[a,b]$, if $s=a$ and $t=b$, and this happens when $\lambda=0$. In particular, if $s=a=1$ and $t=b=-1$, we have the class of Carathoedory functions $P$ consisting of functions $h(z)=1+c_1z+c_2z^2+\cdots$. Many special cases of $P[s,t]$ can also be derived as would be shown later.\vskip 2mm

Since $h=(1+\omega(z))/(1+\omega(z))$, then from (1.2), easy but careful computation shows that:\vskip 2mm

\begin{theorem}
For some $p\in P[s,t]$ and $h\in P$, the following relation holds: $$p(z)=\frac{(1+s)h(z)+1-s}{(1+t)h(z)+1-t}.$$
\end{theorem} 

\begin{theorem}
If $p\in P[s,t]$ then $|p_k|\leq|s-t|$ for $k=1,2,3,\cdots$. Equality holds for $L_0(s,t:z)=(1+sz)/(1+tz)$.
\end{theorem} 

\begin{proof}
Following Goel and Mehrok \cite{RM}, then from (1.2), we have $$1+s\omega(z)=(1+t\omega(z))p(z)=(1+t\omega(z))\left(1+\sum_{k=1}^\infty p_kz^k\right)$$ so that
$$\sum_{k=1}^\infty p_kz^k=\omega(z)\left((s-t)-t\sum_{k=1}^\infty p_kz^k\right).$$  By a lemma of Robertson (see \cite{RM}, Lemma 2), we have
$$\sum_{k=1}^n |p_k|^2\leq|s-t|^2+|t|\sum_{k=1}^{n-1}|p_k|^2.$$ That is,
$$|p_n|^2\leq|s-t|^2-(1-|t|)\sum_{k=1}^{n-1}|p_k|^2\leq|s-t|^2$$ which proves the inequality.
\end{proof}

\begin{theorem}If $p\in P[s,t]$ then $$|p(z)|\leq\frac{1+|s|r}{1-|t|r}.$$\end{theorem}

\begin{proof}
	From (1.2), we have $$p(z)=\frac{1+s\omega(z)}{1+t\omega(z)}.$$ Expanding in series of $\omega(z)$, we have
	$$p(z)=1+(s-t)\omega(z)+t(s-t)\omega(z)^2+t^2(s-t)\omega(z)^3+\cdots.$$ That is,
	$$p(z)=1+(s-t)\sum_{k=1}^\infty t^{k-1}\omega(z)^k$$ so that
	$$|p(z)|\leq 1+|s-t|\sum_{k=1}^\infty |t|^{k-1}|\omega(z)|^k.$$ Since $\omega(z)|\leq|z|=r$, we have
	$$|p(z)|\leq 1+|s-t|\sum_{k=1}^\infty |t|^{k-1}r^k=1+|s-t|\frac{r}{1-|t|r}\leq\frac{1+|s|r}{1-|t|r}$$ as required.
\end{proof}

\begin{theorem}If $p\in P[s,t]$. Define $$q(z)=\int_0^z p(\rho)d\rho.$$ Then $$|q(z)|\leq-\frac{|s|}{|t|}r-\frac{|s|+|t|}{|t|^2}\ln(1-|t|r).$$\end{theorem}

\begin{proof}
Given that $$q(z)=\int_0^z p(\rho)d\rho,$$ then letting $\rho=\sigma e^{i\theta}$. Then we have $$|q(z)|\leq\int_0^\rho |p(\sigma e^{i\theta})|d\sigma$$ Thus we have $$|q(z)|\leq\int_0^{|z|=\rho}\frac{1+|s|\sigma}{1-|t|\sigma}d\sigma$$ which by integration proves the theorem.
\end{proof}

We give examples of special cases of classes of functions packed in $P[s,t]$ for specific choices of $s$ and $t$ defined by (1.1).\vskip 2mm

Following are a few:\vskip 2mm

\noindent 
{\rm(i)} $P[1,-1]\equiv P$, the class of Caratheodory functions for which Re $p(z)>0$,\vskip 2mm

\noindent 
{\rm(ii)} $P[1-2\beta,-1]\equiv P_\beta$, the class of Caratheodory functions of order $\beta$, $0\leq\beta<1$, that is Re $p(z)>\beta$,\vskip 2mm

\noindent 
{\rm(iii)} $P[1,1/\beta-1]\equiv P(\beta)$, the class of functions $p(z)$ satisfying $|p(z)-\beta|<\beta$, $\beta>1/2$,\vskip 2mm

\noindent 
{\rm(iv)} $P[\beta,-\beta]\equiv P^{(\beta)}$, the class of functions $p(z)$ satisfying $|p(z)-1|<\beta|p(z)+1|$, $0<\beta\leq 1$\vskip 2mm

\noindent 
{\rm(v)} $P[\beta,0]\equiv P_{(\beta)}$, the class of functions $p(z)$ satisfying $|p(z)-1|<\beta$, $0<\beta\leq 1$.\vskip 2mm

Furthermore for real number $\lambda\in(-1,1)$ and $\beta$, we have the following examples:\vskip 2mm 

\noindent 
{\rm(vi)} $P\left[\frac{1-i\lambda}{1+i\lambda},\frac{i\lambda-1}{1+i\lambda}\right]$ consisting of functions $p(z)=1+p_1z+\cdots$ such that $$p(z)\prec\frac{1+i\lambda+(1-i\lambda)z}{1+i\lambda+(i\lambda-1)z}.\eqno{(1.3)}$$

\noindent 
{\rm(vii)} $P\left[\frac{1-2\beta-i\lambda}{1+i\lambda},\frac{i\lambda(1-2\beta)-1}{1+i\lambda}\right]$ consisting of functions $p(z)=1+p_1z+\cdots$ such that $$p(z)\prec\frac{1+i\lambda+(1-2\beta-i\lambda)z}{1+i\lambda+(i\lambda(1-2\beta)-1)z},\;0\leq\beta<1.\eqno{(1.4)}$$

\noindent 
{\rm(viii)} $P\left[\frac{\beta+i\lambda(1-\beta)}{\beta(1+i\lambda)},\frac{1-\beta+i\beta\lambda}{\beta(1+i\lambda)}\right]$, $\beta>1/2$, consisting of functions $p(z)=1+p_1z+\cdots$ such that $$p(z)\prec\frac{\beta(1+i\lambda)+[\beta+i\lambda(1-\beta)]z}{\beta(1+i\lambda)+[(1-\beta)+i\beta\lambda]z},\;\beta>1/2.\eqno{(1.5)}$$

\noindent 
{\rm(ix)} $P\left[\frac{\beta(1-i\lambda)}{1+i\lambda},\frac{\beta(i\lambda-1)}{1+i\lambda}\right]$ consisting of functions $p(z)=1+p_1z+\cdots$ such that $$p(z)\prec\frac{1+i\lambda+\beta(1-i\lambda)z}{1+i\lambda+\beta(i\lambda-1)z},\;0<\beta\leq 1.\eqno{(1.6)}$$

\noindent 
{\rm(x)} $P\left[\frac{\beta}{1+i\lambda},\frac{i\beta\lambda}{1+i\lambda}\right]$ consisting of functions $p(z)=1+p_1z+\cdots$ such that $$p(z)\prec\frac{1+i\lambda+\beta z}{1+i\lambda+i\beta\lambda z},\;0<\beta\leq 1.\eqno{(1.7)}$$

\section{New Subclasses of Functions associated with $P[s,t]$}

Let $A$ be the class of functions $f(z)=z+a_2z^2+\cdots$ which are analytic in the unit disk $E$. Denote by $\phi(z)$ the superordinate functions in Equations (1.3) to (1.7), and define $$f'(z)\prec\phi(z),$$
$$\frac{zf'(z)}{f(z)}\prec\phi(z)$$ and
$$1+\frac{zf''(z)}{f'(z)}\prec\phi(z),$$
for functions $f\in A$. Then these classes of functions are respectively new subclasses of functions of bounded turning, starlike and convex functions, which are associated with the extended family $P[s,t]$ of functions having complex parameters.\vskip 2mm

The expositions on new subclasses of starlike and convex functions (and many more function classes) associated with the extended family $P[s,t]$ will be of independent interest. However, in this section, we shall begin the discussion of the first of these new subclasses, that is, new subclasses of functions of bounded turning associated with the extended family $P[s,t]$ of functions having complex parameters, as follows:\vskip 2mm

{\sc Definition.} A function $f\in A$ is said to belong to the class $B(\beta,\lambda)$ if and only if $$f'(z)\prec \frac{1+i\lambda+\beta z}{1+i\lambda+i\beta\lambda z},\;0<\beta\leq 1,\;z\in E.$$

The class $B(\beta,\lambda)$ is well-defined. The following function:\vskip 2mm

$$f(z)=-\frac{i}{\lambda}z+\frac{\lambda^2-i}{\beta\lambda^2}\ln\left(1+\frac{\beta\lambda}{1+i\lambda}z\right)\\
=z+\frac{1}{2}\frac{\beta(i\lambda-1)}{1+i\lambda}z^2+\cdots\eqno{(2.1)}$$ whose derivative is given by $$\phi(z)= \frac{1+i\lambda+\beta z}{1+i\lambda+i\beta\lambda z}$$ is in $B(\beta,\lambda)$.\vskip 2mm

The class $B(\beta,\lambda)$ is a subclass of the class, $R$, of functions whose derivatives have positive real parts in the unit disk $E$ \cite{TH} as we shall prove next:\vskip 2mm

\begin{theorem}
$B(\beta,\lambda)\subset R$, that is, $B(\beta,\lambda)$ consists of univalent functions only.
\end{theorem}

\begin{proof}
It is sufficient to prove that the superordinate function $$\phi(z)= \frac{1+i\lambda+\beta z}{1+i\lambda+i\beta\lambda z},\;0<\beta\leq 1,\;z\in E$$ has positive real parts in $E$, that is, Re $\phi(z)>0$. Let $\alpha(z)=1+\beta z$. Then we have $$2\text{Re}\;\phi(z)=\frac{i\lambda+\alpha}{1+i\lambda\alpha}+\frac{\bar{\alpha}-i\lambda}{1-i\lambda\bar{\alpha}}.$$ Simplifying, we have $$2\text{Re}\;\phi(z)=\frac{2(1+\lambda^2)\text{Re}\;\alpha}{|1+i\lambda\alpha|^2}.$$ Since $z\in E$, it is easy to see that $\text{Re}\;\alpha(z)=\text{Re}(1+\beta z)>0$, and so $\text{Re}\;\phi(z)>0$.
\end{proof}

\begin{theorem}
If $f\in B(\beta,\lambda)$, then $|a_k|\leq\beta/k$. Equality is attained by the function $f(z)$ given by $(2.1)$.
\end{theorem}

\begin{proof}
	Since $f\in B(\beta,\lambda)$, then there exists $p\in P[s,t]$ with $s=\frac{\beta}{1+i\lambda}$ and $t=\frac{i\beta\lambda}{1+i\lambda}$ such that $$f'(z)=p(z).$$ Hence, we have $$|ka_k|\leq|s-t|=\beta$$ so that $|a_k|\leq\beta/k$.
\end{proof}

\begin{theorem}If $f\in B(\beta,\lambda)$ then $$|f'(z)|\leq\frac{1+\lambda^2+\beta r}{1+\lambda^2-\beta\lambda r}.$$ Equality is also attained by the function $f(z)$ given by $(2.1)$.\end{theorem}

\begin{proof}
	Since $f\in B(\beta,\lambda)$, then there exists $p\in P[s,t]$ with $s=\frac{\beta}{1+i\lambda}$ and $t=\frac{i\beta\lambda}{1+i\lambda}$ such that $$f'(z)=p(z).$$ Hence, we have $$|f'(z)|=|p(z)|\leq\frac{1+|s|r}{1-|t|r}$$ so that using $|s|=\beta/(1+\lambda^2)$ and $|t|=\beta\lambda/(1+\lambda^2)$, we have the inequality as desired.
\end{proof}

\begin{theorem} Let $f\in B(\beta,\lambda)$. Then
	\begin{multline*}
	|f(z)|\leq-\frac{1}{\lambda}r-\frac{(1+\lambda)(1+\lambda^2)}{\beta\lambda^2}\ln\left(1-\frac{\beta\lambda}{1+\lambda^2}r\right)\\
	=r+\frac{\beta(1+\lambda)}{2(1+\lambda^2)}r^2+\cdots.
	\end{multline*}
	Equality is also attained by the function $f(z)$ given by $(2.1)$.
\end{theorem}

\begin{proof}
This follows from Theorem 5 taking $f(z)=q(z)=\int_0^z p(\rho)d\rho$ and using $|s|=\beta/(1+\lambda^2)$ and $|t|=\beta\lambda/(1+\lambda^2)$ as required.
\end{proof}

\section{The Case $\lambda=0$}

The class $B(\beta,0)$ is the subclass of functions $f(z)$ such that $$|f'(z)-1|<\beta,\;0<\beta\leq 1,\;z\in E.$$ The following corollaries follow by taking $\lambda=0$ in $B(\beta,\lambda)$:\vskip 2mm

\begin{corollary}
	If $f\in B(\beta,0)$, then $|a_k|\leq\beta/k$. Equality is attained by the function $f(z)$ given by:\vskip 2mm
	$$f(z)=z+\frac{\beta}{2}z^2+\cdots.\eqno{(2.2)}$$
\end{corollary}

\begin{corollary}
	If $f\in B(\beta,0)$ then $$|f'(z)|\leq 1+\beta r.$$
	Equality is also attained by the function $f(z)$ given by $(2.2)$.
\end{corollary}

\begin{corollary}
	Let $f\in B(\beta,0)$. Then $$
	|f(z)|\leq r+\frac{\beta}{2}r^2+\cdots.$$
	Equality is also attained by the function $f(z)$ given by $(2.2)$.
\end{corollary}

\section{Concluding Remarks}

The study of many subclasses of functions opened up in this paper such as the subclasses of starlike and convex functions associated with Janowski functions $P[s,t]$ having complex parameters shall make fruitful engagements for interested researchers in the near and distant future.\vskip 2mm

{\it Acknowledgements.} This work is in honour of my teacher, supervisor, mentor and professor of Mathematics, Professor T. O. Opoola, who retired recently from active service from the University of Ilorin, Ilorin, Nigeria. The author appreciates the professor's immense role in his academic development, especially believing in his ability despite years in the industry before turning back to begin a career in the academia.\vskip 2mm

This work was carried out at the Department of Mathematics, Obafemi Awolowo University, Ile-Ife, Nigeria during the author's visiting position in the Department. The author is indebted to all staff of the Department for their hospitality, especially the head, Professor B. S. Ogundare.\vskip 2mm

The author is indebted to the referee for helpful comments and suggestions.

\vspace{10pt}

\hspace{-4mm}{\small{}}

\vspace{-12pt}
\begin{tabular}{l}
{\small  {\bf Current Address}}\\
{\small  Department of Mathematics}\\
{\small  Obafemi Awolowo University, Ile-Ife, Nigeria}\\
{\small {\tt kobabalola@oauife.edu.ng}}\\
\end{tabular}

\vskip 2mm

\begin{tabular}{l}
	{\small  {\bf Permanent Address}}\\
	{\small  Department of Mathematics}\\
	{\small  University of Ilorin, Nigeria}\\
	{\small {\tt babalola.ko@unilorin.edu.ng}}\\
	{\small {\tt kobabalola@gmail.com}}\\
\end{tabular}

\end{document}